# MORE RIGOROUS RESULTS ON THE KAUFFMAN–LEVIN MODEL OF EVOLUTION

By Vlada Limic[1] and Robin Pemantle[2]

*University of British Columbia and University of Pennsylvania*

The purpose of this note is to provide proofs for some facts about the NK model of evolution proposed by Kauffman and Levin. In the case of normally distributed fitness summands, some of these facts have been previously conjectured and heuristics given. In particular, we provide rigorous asymptotic estimates for the number of local fitness maxima in the case when $K$ is unbounded. We also examine the role of the individual fitness distribution and find the model to be quite robust with respect to this.

**1. Introduction.** The purpose of this note is to provide proofs for some facts about the NK model. Some of these proofs have been previously formulated, at least approximately, as conjectures or heuristic arguments. Since we are interested in the mathematical analysis of the model, we include only a brief summary of the biological motivation, for which we can do no better than to excerpt and paraphrase from the introductory section of the paper by Evans and Steinsaltz [3].

Beginning with Sewall Wright in the early twentieth century, evolution has been modeled as the gradual motion of a genome through an abstract space, with a tendency toward increasing values of the *fitness function*. One may think of the graph of this function as a *fitness landscape* and of natural selection as a random walk with upward drift on the fitness landscape. One cannot understand the likely behavior of such a random walk without understanding the qualititative nature of the landscape as one with "slivers of high fitness looming up above the vast genomic tohubohu" [3]. In any random walks model of fitness landscapes and natural selection, the nature of the global fitness maximum is less important than the number and height of local maxima.

Received November 2002; revised August 2003.

[1]Supported in part by NSF Grant DMS-01-04232 and by NSERC research grant.

[2]Supported in part by NSF Grant DMS-01-03635.

Kauffman and Levin [7] introduced the NK model, which is a probabilistic model for the fitness landscape. In this model, there are $N$ loci, at each of which is one of two possible alleles. Thus a genome is an element of the space $\{0, 1\}^N$. The fitness of a genome is the sum of $N$ different fitnesses, the $j$th of which is determined by the alleles at sites $j, j + 1, \ldots, j + K$ modulo $N$. In the NK model, the $2^{K+1}$ alleles in the $N$ possible positions are given fitnesses whose joint distribution is that of $2^{K+1}N$ i.i.d. picks from a distribution $F$. The fitness of a given genome is then the sum of the $N$ fitnesses corresponding to the actual string of $K + 1$ alleles beginning at each position. Note that this randomness is present in the model at the start; later one may model natural selection as a random walk in this random environment, but that is beyond the scope of this paper. Evans and Steinsaltz pointed out that since the allele substrings of length $K + 1$ overlap, there is no easy way to find the optimal choice for the $N$ alleles. They concluded that "while no one would mistake this abstract system for a realistic model of genetic evolution, it has the virtues of a good foundational model: it is easy to describe, yet contains a wealth of structure that is neither obvious nor superficially accessible. Before we can analyze a more realistic model, it would seem we must first come to grips with models such as this one. At the same time, we may hope that some general features of this model will carry over to something like the real world."

Most studies of the NK model rely on simulations, which are limited to small to intermediate values of $N$ (e.g., in  [6], $N = 96$ and in [2], $N = 1024$, which corresponds to the size of a gene, but it is much smaller than the number of genes in a genome). Simulations may provide quick answers to various questions in particular cases of fitness distribution $F$. However, a very interesting and natural question of robustness of the model under variations in $F$ can be tackled only mathematically.

We warn the reader that we *always* assume in this paper that the parameter $K$ is strictly positive and that the underlying distribution $F$ is continuous. The NK model for $K = 0$ or $K = N - 1$ exhibits special behaviors which were rigorously analyzed by many authors (see, e.g., [7]). If $F$ were not continuous, ties would be possible and analysis would become cumbersome.

The study of the question to which our paper is devoted begins with [11], where Weinberger gives asymptotic formulae for the number of local fitness maxima (LFM) when $N$ and $K$ are large and $F$ is the normal distribution. As noted in [2] and [3], however, Weinberger's derivation is not rigorous. Weinberger's heuristics are limited to the case where $F$ is the normal distribution, although he points out that other distributions such as the Cauchy might be more realistic and that one could expect the outcome to be independent of the choice of distribution.

The majority of rigorous results that have been obtained assume that $K$ is fixed and $N \to \infty$. In this context, several results were obtained in two recent papers [2, 3]. Among other things, they both show ([3], Theorem 7, and [2], Theorem 2.1) that the exponential growth rate number of local maxima (or, equivalently, the exponential decay of the probability of a given genome being a local fitness maximum) exists as a limit. In other words, the probability of a LFM decays like $\exp N(\lambda_K + o(1))$ as $N \to \infty$ with $K$ remaining fixed. For $K = 1$, they computed this limit explicitly when $F$ is the exponential distribution [3] or the negative exponential [2]. In the case where $F$ has an exponential moment, Durrett and Limic ([2], Theorem 5.1) made partial progress toward showing the number of local maxima (for large $K, N$) to be independent of the distribution $F$: they bounded the exponential rate on one side and they conjectured this to be correct to within a constant factor. The value of $\lambda_K$ is theoretically possible to compute for certain distributions when $K \geq 1$, but practically impossible. It is biologically reasonable that $K$ be on the order of at least several dozens, whence our interest in asymptotic formulae for $\lambda_K$ with error estimates that are valid as $K, N \to \infty$ without restriction. For example, in [6], pages 122–142, it is shown that maturation of the immune response fits the parameters $K = 40$ and $N = 122$, which is probably best described as "$N$ and $K$ large, with $N/K$ remaining bounded."

The first purpose of this note is to rigorize Weinberger's computations for the normal case. This includes sharpening his statements to include error bounds and quantified asymptotic statements, specifically convergence uniform in $N$ as $K \to \infty$. The second purpose is to investigate dependence on $F$. Specifically, we prove some asymptotic results that do not depend at all on the distribution of $F$, completing and generalizing the conjecture in [2], and we show some stronger results for the "fat-tail" case, which we believe to be the extreme opposite to the case where $F$ has finite second moment.

The remainder of the paper is organized as follows. The next section sets forth the notation and states our main results. Section 3 gives proofs for the results in which $F$ is the normal distribution. Section 4 proves results for general distributions and derives asymptotics for fat-tailed distributions when $N/K \to \infty$. Section 5 contains a detailed analysis of the case where $F$ has fat tails and $N/K$ remains bounded. Finally, Section 6 gives an exact expression for the exponential rate when $F$ is the fat tail and $K = 1$, which, when compared with similar computations for other distributions, corroborates an extremality conjecture for the fat tail.

We use notation $o(1)$ to represent a term that converges to 0 as $K \to \infty$, $O(1)$ to represent a term bounded by a constant and $\Theta(\text{expression}(K))$ to represent a term for which there are positive finite constants $c, C$ (independent of $K$) such that $c \, \text{expression}(K) \leq \text{term} \leq C \, \text{expression}(K)$.

**2. Notation and statements of results.** The parameters of the model are positive integers $N > K$ and a continuous distribution function $F$ on the real numbers. Our concern in this paper is with the number of LFMs for a random fitness landscape. The expectation of this number is equal to $2^N$ times the probability that any given genome is a local fitness maximum. Consequently, our sole focus is the rigorous estimation of this probability. Showing that the logarithm of the number of LFMs is near its expectation is not hard, but will not concern us here; see, for example, [2], Theorem 7.1, where an asymptotic normality result is obtained for the logarithm of the number of local fitness maxima.

In the NK model the (unnormalized) fitness of a particular genome $\eta = (\eta_1, \eta_2, \ldots, \eta_N) \in \{0, 1\}^N$ is defined to be

$$\sum_{j=1}^N Y(j; (\eta_j, \eta_{j+1}, \ldots, \eta_{j+K})), \tag{2.1}$$

where the family

$$\{Y(j; (\eta_1, \eta_2, \ldots, \eta_{K+1})) : j = 1, \ldots, N; (\eta_1, \eta_2, \ldots, \eta_{K+1}) \in \{0, 1\}^{K+1}\}$$

is the family of of $N \cdot 2^{K+1}$ i.i.d. random variables with common distribution $F$. Suppose we are given such a family on a probability space $(\Omega, \mathcal{F}, \mathbb{P})$ and abbreviate

$$Y_j := Y(j; (0, 0, \ldots, 0))$$

to be the fitness of the substring of $K + 1$ zeros starting in position $j$; here and *throughout*, arithmetic on subscripts is always taken modulo $N$. With the above notation the fitness of the zero genome is $\sum_{j=1}^N Y_j$.

The genome consisting of all 0's has $N$ neighbors, namely all binary strings of length $N$ with exactly one 1. Since in this paper we are only interested in the probability of the event that the string of all 0's is LFM, the only other relevant random variables from the above family are the fitnesses $Y(j; (\eta_1, \eta_2, \ldots, \eta_{K+1}))$, where $j = 1, \ldots, N$ and where $\sum_i \eta_i = 1$. We again abbreviate for $1 \leq j \leq N, 0 \leq i \leq K$,

$$Y_{j,i} := Y(j - i; (0, \ldots, 1, \ldots, 0)),$$

where 1 is only in the $i$th position above (here we count positions starting from 0). The quantity $Y_{j,i}$ is interpreted as the fitness of the substring of length $K + 1$ starting at position $j - i$ that is all 0's except for a single 1 in position $j$. Then the definition (2.1) says that the string $\mathbf{e}_j$ consisting of $N - 1$ 0's and a single 1 in the $j$th position has fitness (in the new notation)

$$\sum_{i=j+1}^{j-K-1} Y_i + Y_{j,0} + Y_{j,1} + \cdots + Y_{j,K}.$$

The zero genome is a LFM if it has greater fitness than that of any genome with exactly one 1. We denote the event of optimality of the zero string by $\mathcal{H}$. We may write $\mathcal{H} = \bigcap_j \mathcal{H}_j$, where $\mathcal{H}_j$ is the event that all 0's are better than $\mathbf{e}_j$. Equivalently,

$$(2.2) \qquad \mathcal{H}_j \Leftrightarrow \sum_{i=j-K}^{j} Y_i \geq \sum_{i=0}^{K} Y_{j,i}.$$

Define $p_F(N, K) := \mathbb{P}(\mathcal{H})$. We usually suppress dependence on $F$ and write simply $p(N, K)$. Our first result makes rigorous and precise what is stated in [11].

THEOREM 2.1. *Suppose that $F$ is the standard normal distribution. Then*

$$\log p(N, K) = \frac{N}{K}(-\log K + R_{N,K})$$

*with*

$$c \log \log K \geq R_{N,K} \geq -c\sqrt{\log K}$$

*for some $c > 0$.*

REMARKS. (i) Specializing to the case $N/K \to \alpha$, we obtain the estimate $p(N, K) = K^{-1/\alpha + o(1)}$. (ii) The error terms are independent of $N$, so the previous estimate is uniform in $N > K + 1$ as $K \to \infty$; here and throughout, *all asymptotic notation is with respect to $K$ only* (unless otherwise noted). (iii) In contrast to what will be the case with other distributions, there is no correction when $N/K$ does not go to infinity. (iv) If $K = N - 1$, then the NK model is essentially different from the NK model where $K < N-1$, but since $p(N, N-1) = 1/N + 1$, it is still true that $\log p(N, N-1) = -\log(N+1) \sim -N \log(N-1)/(N-1)$ with error smaller than the above bounds on $R_{N,N-1}$ for large $N$.

Next, we state our most general result.

THEOREM 2.2. *Let $F$ be any distribution and $N \geq 2(K+1)$. Then*

$$(2.3) \qquad \log p(N, K) \leq -(1 + o(1)) \left\lfloor \frac{N}{K} - o(1) \right\rfloor \log K$$

$$(2.4) \qquad \geq -(3 + o(1)) \left\lceil \frac{N}{K} \right\rceil \log K.$$

We believe that the upper bound (2.3) is sharp, so we make the following conjecture:

CONJECTURE 1. *It is possible to replace* 3 *by* 1 *in* (2.4).

When sums of random variables are concerned, the class of most tightly clustered distributions comprises the distributions with finite variance, since these exhibit Gaussian behavior when summed. At the other extreme, one has distributions with extremely fat tails. In the limit, one might consider a distribution with the following property: In any collection of $n$ i.i.d. picks, the greatest is much greater than the sum of the magnitudes of the others with probability tending exponentially rapidly to 1 as $n \to \infty$. For example, if $U$ is uniform on $[0, 1]$, then $\exp(\exp(1/U))$ has this property. In this case, as long as $K \to \infty$ at least as fast as $\log N$, one may approximate $\mathcal{H}_j$ by the event

$$(2.5) \qquad \mathcal{H}'_j := \left\{ \max_{j-K \leq i \leq j} Y_i \geq \max_i Y_{j,i} \right\}.$$

Heuristically, properties of $p(N, K)$ shared by fat-tailed distributions and normal distributions would be likely to hold for all distributions, since all others lie in between. One approach to establishing facts about fat-tailed distributions would be to axiomatize how fast the probability should tend to 1 of the event that the largest of $n$ picks dominates all the others, and then prove theorems about distributions satisfying the axiom. We choose a less cumbersome approach, namely to provide an analysis of the probability of the event $\mathcal{H}' := \bigcap_{j=1}^{N} \mathcal{H}'_j$. We use the notation $p_{\text{fat}}(N, K)$ to denote $\mathbb{P}(\mathcal{H}')$ and sometimes call it "$p(N, K)$ under the fat-tail distribution." Note that $p_{\text{fat}}(N, K)$ is independent of $F$, assuming $F$ is continuous.

CONJECTURE 2. *For any* $N$ *and* $K$, *the infimum over all* $F$ *of* $p_F(N, K)$ *is equal to* $p_{\text{fat}}(N, K)$.

Our next result shows that Conjecture 1 holds for the fat tail and thus that Conjecture 2 implies Conjecture 1.

THEOREM 2.3. *We have*

$$\log p_{\text{fat}}(N, K) \geq -(1 + o(1)) \left\lceil \frac{N}{K} + o(1) \right\rceil \log K.$$

Weinberger suggested the Cauchy as a biologically realistic distribution. Those readers who are bothered by a mythological distribution called the fat tail will perhaps be interested to see that the previous result for the fat tail may be proved for the Cauchy. We remark that the criterion we have suggested for axiomatization of the fat tail, namely exponential decay of the probability that the largest of $n$ picks fails to dominate the sum of the others, requires much fatter tails than the Cauchy distribution possesses. Thus we view the following result as more than adequate to demonstrate that the fat-tail results hold for typical fat-tailed distributions.

THEOREM 2.4. *When $F$ is a symmetric Cauchy distribution,*

$$\log p(N, K) \geq -(1 + o(1)) \left\lceil \frac{N}{K} + o(1) \right\rceil \log K.$$

Comparing these last results to Theorem 2.1, we see that for the fat-tail and Cauchy distributions, and conjecturally for all distributions $F$, $\log p_F(N, K) \sim \log p_\Phi(N, K)$, where $\Phi$ is the normal c.d.f., as long as $N/K \to \infty$: In this case the difference between $N/K$ and $\lceil N/K + o(1) \rceil$ is irrelevant and the formulae agree. Note that, on the other hand, if $N/K \approx \alpha$, where $\alpha = m - 0.5$ for some integer $m$, the difference between $\lceil N/K + o(1) \rceil$ and $\lfloor N/K + o(1) \rfloor$ is 1, which amounts to the difference of $1/K$ in the asymptotic lower and upper bounds for $p(N, K)$. It turns out there is, in fact, an asymptotic inequivalence between $\log p_\Phi(N, K)$ and $\log p_{\text{fat}}(N, K)$ when $N/K$ does not go to infinity. Because of this, we include a more precise description of that asymptotics of $\log p(N, K)$ in this regime.

The statement of the following theorem makes more sense if one keeps in mind how $\mathcal{H}'$ is likely to occur. There will be at least $r_0 := \lceil N/K \rceil$ large fitnesses among the $Y_j$, which is the minimum number for which it is possible to have a large fitness in every window of size $K$. The number of ways to pick $r$ large fitnesses increases with $r$, but the probability that any specific $r$ fitness values are all large decreases with $r$. In this energy–entropy tradeoff, the maximum occurs at $r = r_0$ as $N/K$ increases to $r_0 - o(1)$, at which point the $r$ value that achieves the maximum switches to $r_0 + 1$.

THEOREM 2.5. *As $K \to \infty$ with $N/K$ bounded, there are formulae that give the value of $p_{\text{fat}}(N, K)$ up to a factor of $1 + o(1)$. The formulae are in terms of functions $\{f_r : r \geq 3\}$ on $\mathbb{R}^+$, which are defined by formula (5.7) in Section 5 and summarized in Table 1. Additionally, the functions $f_r$ satisfy the following statements:*

- $f_r(0) = 0$.
- $f_r(x) \sim x^{r-1}$ *as $x \to 0$.*
- *For $r \geq 4$, $f_r$ is increasing, continuous and bounded on $[0, 1]$.*
- *For $r = 3$, $f_3$ is increasing and continuous on $[0, 1)$, with $f_3(1 - t) \sim 2 \log(1/t)$ as $t \to 0^+$.*

In other words, there are narrow windows in the parameter $N/K$ in which $p_{\text{fat}}(N, K)$ changes from roughly $K^{-r}$ to $K^{-(r+1)}$. These windows occur at $N/K \approx r - K^{-1/(r-1)}$. An exception is when $r = 2$. In this case, the change from order $K^{-2}$ to order $K^{-3} \log K$ is complete at $N = 2K - c \log K$, after which the order slowly slides down to $K^{-3}$ as $\log(N - 2K)$ increases to $\log K$.

TABLE 1
*Behavior of $p(N, K)$ across integer values of $N/K$*

| $N$ | $j$ | $(1 + o(1))p_{\mathrm{fat}}(N, K)$ |
|---|---|---|
| $2(K+1) - j$ | $0 < j \le K$ | $\frac{1}{K}\left(\frac{1}{K+3-j} - \frac{1}{K} + \frac{\log K}{K^2}\right)$ |
| $(r - y)(K + 1)$ | $0 \le y < 1$ | $\frac{1}{K^r}f_r(y) + \frac{1}{K^{r+1}}f_{r+1}(1 + y)$ |

A final result is the analysis for the fat tail when $K = 1$. Note that when $K = O(1)$, maxima are taken over collections of a bounded size, so no actual distribution has tails fat enough to ensure that the maximum dwarfs the others. Nevertheless, this result is still relevant to Conjecture 2.

THEOREM 2.6. *We have*

$$N^{-1}\log p_{\mathrm{fat}}(N, 1) \to z := -\log 1.803\ldots = -0.58947\ldots,$$

*where $z$ is the solution of the Bessel equation*

$$0 = \pi\sqrt{6}\,\mathrm{BesselI}(\tfrac{2}{3}, \tfrac{2}{3}\sqrt{2z}) - \pi\sqrt{3z}\,\mathrm{BesselI}(-\tfrac{1}{3}, \tfrac{2}{3}\sqrt{2z})$$
$$+ 3\sqrt{2}\,\mathrm{BesselK}(\tfrac{2}{3}, \tfrac{2}{3}\sqrt{2z}) + 3\sqrt{z}\,\mathrm{BesselK}(\tfrac{1}{3}, \tfrac{2}{3}\sqrt{2z}).$$

The published exact values of $\log p(N, 1)$ for the exponential and negative exponential are, respectively, $-0.57504\ldots$ [3] and $-0.5499934\ldots$ [2]. The published lower bound for the uniform is $-0.55957\ldots$ [2]. All of these values are greater than the value for the fat tail given by Theorem 2.6, thus providing further corroboration of Conjecture 2.

Some final notation and methodology common to all the proofs is as follows. We let $\mathcal{F} = \sigma(Y_j : 1 \le j \le N)$ be the $\sigma$-field generated by the fitnesses of zero substrings. We let $F^{(K+1)}$ denote the c.d.f. for the sum of $K + 1$ independent picks from the distribution $F$. Conditional on $\mathcal{F}$, the events $\mathcal{H}_j$ are independent, with

$$\mathbb{P}(\mathcal{H}_j | \mathcal{F}) = F^{(K+1)}\left(\sum_{i=j}^{j+K} Y_i\right).$$

Removing the conditioning then gives a formula which appears as [11], (2.4),

$$(2.6) \qquad p(N, K) = \int \prod_{j=1}^{N} F^{(K+1)}\left(\sum_{i=j}^{j+K} Y_i\right) dF(Y_1) \cdots dF(Y_N).$$

**3. Analysis of the normal case.** The following facts are well known.

LEMMA 3.1. *If $\Phi$ and $\phi$ are the normal c.d.f. and density, respectively, then*

(3.1)
$$\log \Phi(x) = (-1 + o(1))(1 - \Phi(x))$$
$$= \phi(x)(x^{-1} + O(x^{-2})), \qquad x \to \infty,$$

(3.2)
$$(\log \Phi)'' = \frac{\Phi \phi' - \phi^2}{\Phi^2} < 0$$

*and*

(3.3)
$$\textit{the function } \log \Phi \textit{ is concave.}$$

Next we define the normalized total fitness

$$t := N^{-1/2} \sum_{j=1}^{N} Y_j$$

and the recentered window sums

$$X_j := \frac{(\sum_{i=j}^{j+K} Y_j) - ((K+1)/\sqrt{N})t}{\sqrt{K+1}}.$$

It is immediate to verify that each $X_j$ is a normal with mean 0 and variance $1 - (K + 1/N)$. Since the quantities $Y_j - t/\sqrt{N}$ are independent normals recentered to sum to zero, their joint distribution is independent of the centering constant $t$. This can be verified explicitly by checking that the covariance of $X$ and $Y_j - t/\sqrt{N}$ is 0 for each $j$. Consequently, since $\sqrt{K+1}X_j = \sum_{i=j}^{j+K}(Y_j - t/\sqrt{N})$, we see that

(3.4)
$$\{X_j : 1 \leq j \leq N\} \text{ is independent of } t.$$

Plugging this into (2.6) and using the fact that $F^{(K+1)}$ is a normal of variance $K + 1$, we get

(3.5)
$$p(N, K) = \mathbb{E} \prod_{j=1}^{N} \Phi\left(X_j + \sqrt{\frac{K+1}{N}}t\right).$$

Up to here we have followed Weinberger, arriving at [11], (3.2). Weinberger now asserts that $X_j = O(1)$ with mean zero, and may therefore be removed from the equation, resulting in $p(N, K) \approx \mathbb{E}\Phi(t\sqrt{(K+1)/N})^N$, where $t$ is a standard normal; this is then evaluated by steepest descent. Our contribution in the rest of this section is to finish this properly, with one inequality (the upper bound on $R$) following directly from (3.3) of Lemma 3.1, rather than relying on independence of $t$ and $\{X_j : 1 \leq j \leq N\}$.

*Upper bound on* $R$. By definition, the random variables $X_j$ sum to zero. Using concavity of $\log \Phi$, we have the (deterministic) inequality

$$\sum_{j=1}^{N} \log \Phi\left(X_j + t\sqrt{\frac{K+1}{N}}\right) \leq N \log \Phi\left(t\sqrt{\frac{K+1}{N}}\right).$$

Plugging into (3.5) then gives

$$(3.6) \qquad \mathbb{P}(A) \leq \mathbb{E}\Phi\left(t\sqrt{\frac{K+1}{N}}\right)^N = \int \Phi\left(x\sqrt{\frac{K+1}{N}}\right)^N \phi(x)\,dx,$$

where $\phi$ is the normal density. Let $I(x) = I_{N,K}(x)$ denote the integrand in (3.6) and let $M$ denote the maximum value of $\log I$:

$$M := \max_x \log I(x) = -\log\sqrt{2\pi} + \max_x \left[N \log \Phi\left(x\sqrt{\frac{K+1}{N}}\right) - \frac{x^2}{2}\right].$$

If we can show that

$$(3.7) \qquad \log \int I_{N,K}(x)\,dx \leq M + O(1)$$

and that

$$(3.8) \qquad M = -\frac{N}{K}(\log K + O(\log\log K)),$$

then the first inequality in Theorem 2.1 will be proved. Both computations are routine, and we need only one inequality of (3.8), but we include the arguments because they clarify matters by indicating the location of the saddle.

To show (3.8), let $x_0 := \sqrt{(2N/(K+1))\log(K+1)}$. Of course

$$M \geq \log I_{N,K}(x_0)$$

$$= -\log\sqrt{2\pi} + \frac{N}{K+1}[(K+1)\log\Phi(\sqrt{2\log(K+1)}) - \log(K+1)]$$

$$= -\frac{N}{K}\left[\log K + \frac{1+o(1)}{\sqrt{\log K} + o(1)}\right],$$

where we have used the estimate (3.1) from Lemma 3.1 on $\log \Phi$ and where the last $o(1)$ accounts for $-\log\sqrt{2\pi}$. This shows one inequality in (3.8). For an upper bound on $M$, suppose first that $x \geq \sqrt{(2N/(K+1))} \times \sqrt{(\log(K+1) - 2\log\log(K+1))}$. Then

$$\log I_{N,K}(x) \leq -\frac{x^2}{2} = -\frac{N}{K}(\log K + O(\log\log K))$$

as needed. On the other hand, when $x \leq \sqrt{(2N/(K+1))} \times \sqrt{(\log(K+1) - 2\log\log(K+1))}$, then

$$\log I_{N,K}(x) \leq -\log\sqrt{2\pi} + N\log\Phi(\sqrt{2(\log(K+1) - 2\log\log(K+1))}\,)$$

$$= (-1 + o(1))N\frac{1}{K}\frac{(\log K)^2}{\sqrt{2\log K - 4\log\log K}}$$

$$\leq -(1 + o(1))\frac{N}{K}(\log K)^{3/2},$$

so these values of $x$ need not be considered and the other inequality in (3.8) is proved.

Proving (3.7) is merely a matter of estimating the second derivative of $\log I$. By log concavity of $\Phi$, this is at most the second derivative of $\log\phi$, which is equal to $-1/2$. Let $x_M := x_M(N, K)$ be such that $I_{N,K}(x_M) = M$. Now an easy calculus argument (using log concavity) shows

$$I(x) \leq \exp\{I(x_M)\}\exp\{-(x - x_M)^2/4\} = e^M\exp\{-(x - x_M)^2/4\},$$

so that $\int e^{-M}I_{N,K}(x)\,dx$ is bounded above by a constant $2\int_0^\infty \exp(-x^2/4)\,dx$ that is independent of $N$ and $K$, which shows that $\log\int I(x)\,dx \leq M + O(1)$ and finishes the proof of (3.7) and the first inequality of Theorem 2.1.

*Lower bound on $R$.* Let $G_1$ be the event that

$$t \geq x_1 := \sqrt{(2N/(K+1))(\log(K+1) + 3\sqrt{\log(K+1)})}.$$

Let $G_2$ be the event that $\max|X_j| \leq 1$. Due to independence of $t$ from $\{X_j : 1 \leq j \leq N\}$, we may write

$$p(N, K) \geq \mathbb{P}(G_1 \cap G_2)\mathbb{P}(\mathcal{H}|G_1, G_2) = \mathbb{P}(G_1)\mathbb{P}(G_2)\mathbb{P}(\mathcal{H}|G_1, G_2).$$

We estimate this in pieces, the first being the one responsible for pushing $R$ down to $-c\sqrt{\log K}$.

Since $\log\phi(x) = -x^2/2 + O(1)$, we may estimate

$$\log\mathbb{P}(G_1) = \log(1 - \Phi(\sqrt{(2N/(K+1))(\log(K+1) + 3\sqrt{\log(K+1)})}\,))$$

$$= \log\left((1 + o(1))\frac{\phi}{x}[\sqrt{(2N/(K+1))(\log(K+1) + 3\sqrt{\log(K+1)})}\,]\right)$$

$$= O(1) - \frac{N}{(K+1)}(\log(K+1) + 3\sqrt{\log(K+1)}\,) - \log x_1$$

$$= -\frac{N}{K}(\log K + O(\sqrt{\log K})\,).$$

Next, we estimate $\mathbb{P}(G_2)$.

Lemma 3.2. *We have*

$$\log \mathbb{P}(G_2) \geq -\frac{\pi^2}{2}\frac{N}{K}.$$

Proof. Let $S_j := \sum_{i=1}^{j}(Y_i - t/\sqrt{N})$ be the recentered partial sums. Then $X_j = (K+1)^{-1/2}(S_{j+K} - S_{j-1})$, with indices still taken modulo $N$. The event $G_2'$, defined by

$$G_2' := \{|S_j| \leq \tfrac{1}{2}\sqrt{K} \text{ for all } j \leq N\},$$

implies the event $G_2$. Let $W_0$ be Wiener measure on continuous paths $\omega$ on $[0, N]$ starting at 0 and let $W_0^{\mathrm{br}}$ be the Brownian bridge measure, that is, $W_0$ conditioned on $\{\omega(N) = 0\}$. The law of $\{S_j : 1 \leq j \leq N\}$ is the law of partial sums of $N$ i.i.d. standard normals conditioned on summing to zero; this is the same as the conditional law of $\{\omega(j) : 1 \leq j \leq N\}$ under $W_0$, conditioned on $\{\omega(N) = 0\}$, which is the same as the law of $\{\omega(j) : 1 \leq j \leq N\}$ under $W_0^{\mathrm{br}}$.

A Brownian bridge always stays closer to the origin than unconstrained Brownian motion, in the following sense. In fact, it is not difficult to couple the path of the reflected simple random walk bridge (i.e., the absolute value of the random walk path conditioned to visit 0 at time $2n$) and the path of the reflected simple random walk up to step $2n$ so that the former stays below the later at all times with probability 1. Taking the diffusion limits in an appropriate way constructs one coupling of the reflected Brownian bridge and reflected Brownian motion described above.

Letting $G_2''$ be the event that $|\omega(t)| \leq \sqrt{K}/2$ for all $t \leq N$, we then have

$$\mathbb{P}(G_2) \geq \mathbb{P}(G_2') \geq W_0^{\mathrm{br}}(G_2'') \geq W_0(G_2'').$$

Let $W_\mu$ be Wiener measure started from distribution $\mu$. Clearly $W_\mu(G_2'')$ is maximized when $\mu = \delta_0$; that is, $W_0(G_2'') \geq W_\mu(G_2'')$ for any $\mu$. Now let $\mu$ be the distribution on $[-\sqrt{K}/2, \sqrt{K}/2]$ with density $C\cos(\pi x/\sqrt{K})$. This is an eigendensity for Brownian motion killed on exiting $[-\sqrt{K}/2, \sqrt{K}/2]$ (see [8], Theorem 4.1.1). We see that

$$\mathbb{P}(G_2) \geq W_\mu(G_2'') = \exp\left(-\frac{\pi^2}{2}\frac{N}{K}\right),$$

proving the lemma. □

Finally, we estimate the third term. Recall from (3.5) the formula for the probability of LFM:

$$p(N, K) = \mathbb{E}\left[\exp\left(\sum_{j=1}^{N} \log \Phi(X_j + \sqrt{(k/N)}t)\right)\right].$$

For $x > 1$, consider the inequality

$$\sqrt{2(x + 3\sqrt{x}\,)} > \sqrt{2(\sqrt{x} + 1)^2} - \sqrt{2} + 1 = \sqrt{2x} + 1,$$

which can easily be checked, for example, by squaring both sides (note that if $x > 1$, then both sides of the inequality are strictly positive). Applying this inequality yields on $G_1 \cap G_2$,

$$X_j + t\sqrt{\frac{K+1}{N}} \geq \sqrt{2(\log(K+1) + \sqrt{\log(K+1)}\,)} - 1 \geq \sqrt{2\log(K+1)}.$$

Therefore, on $G_1 \cap G_2$ we then have for all $j$

$$\log \Phi\left(X_j + t\sqrt{\frac{K+1}{N}}\right) \geq \log \Phi(\sqrt{2\log(K+1)}\,)$$

and hence, using (3.5),

$$\mathbb{P}(\mathcal{H}|G_1, G_2) \geq (\Phi(\sqrt{2\log(K+1)}\,))^N \geq \exp\left(-\frac{N}{K+1}\right).$$

Plugging in the estimates for $\mathbb{P}(G_1)$ and $\mathbb{P}(G_2)$ then yields

$$\log p(N, K) \geq -\frac{N}{K}\left(\frac{\pi^2}{2} + \log K + O(\sqrt{\log K})\right),$$

which finishes the proof of the theorem.

## 4. Proof of universality results.

PROOF OF THEOREM 2.2. *First inequality.* For the moment let the small positive real parameter $y$ be unspecified. Break the interval from 1 to $N$ into $L := \lfloor N/(1+y)(K+1) \rfloor$ intervals of length $\lfloor (1+y)(K+1) \rfloor$, discarding any unused positions at the end. Denote these intervals $I_1, \dots, I_L$ and let $I'_j$ denote the first $\lceil y(K+1) \rceil$ positions in $I_j$. Let $s_j$ denote the index $s \in I'_j$ that maximizes $S := \sum_{i=0}^{K} Y_{s,i}$. The maximum is a maximum of $y(K+1)$ independent draws from $F^{(K+1)}$, so $B_j := F^{(K+1)}(\sum_{i=0}^{K} Y_{s_j, i})$ has distribution $\beta(1, y(K+1))$. The mean of $B_j$ is $1 - (y(K+1) + 1)^{-1}$. For the event $\mathcal{H}_{s_j}$ to occur, the sum $\sum_{i=s_j}^{s_j+K} Y_i$ must exceed $S$. Let $\mathcal{F}' = \sigma(Y_{j,i} : 1 \leq j \leq N, 0 \leq i \leq K)$ be the $\sigma$-field generated by the fitnesses of substrings with exactly one 1. Then

$$P(\mathcal{H}_{s_j}|\mathcal{F}') = 1 - B_j.$$

Since $|s_j - s_k| > K$ when $j \neq k$, the events $\mathcal{H}_{s_j}$ are conditionally independent given $\mathcal{F}'$, and the $B_j$'s are mutually independent random variables.

Therefore,

$$\mathbb{P}(\mathcal{H}) \leq \mathbb{P}\left(\bigcap_{j=1}^{L} \mathcal{H}_{s_j}\right) = \mathbb{E}\mathbb{P}\left(\bigcap_{j=1}^{L} \mathcal{H}_{s_j} | \mathcal{F}'\right) = \mathbb{E}\prod_{j=1}^{L}(1 - B_j) = \left(\frac{1}{1 + y(K+1)}\right)^L.$$

When $K = o(N)$, we choose $y = y(K) = o(1)$ to optimize this bound. For example, taking $y = 1/\log K$ gives an upper bound of $\exp(-(1 + o(1)) \times (N/K)\log K)$, as is required to prove (2.3).

When $K = \Theta(N)$, the same choice of $y$ leads to the same conclusion, except that one has $\lfloor N/y(K+1) \rfloor$ in place of $N/K$. Since $y(K) = o(1)$, this is again sufficient to prove (2.3).

*Second inequality.* To prove (2.4), begin with the observation that the events $\mathcal{H}_j$ are increasing events with respect to the variables $\{Y_j : 1 \leq j \leq N\}$ and $\{-Y_{j,i} : 1 \leq j \leq N, 0 \leq i \leq K\}$. By Harris' inequality, these are positively associated. Let $L = \lceil N/(K+1) \rceil$ and, for $1 \leq j \leq L$, let

$$G_j := \bigcap_{i=(j-1)(K+1)+1}^{j(K+1)} \mathcal{H}_i.$$

Positive association implies that

$$\mathbb{P}(\mathcal{H}) = \mathbb{P}\left(\bigcap_{j=1}^{L} G_j\right) \geq \mathbb{P}(G_1)^L.$$

Thus it suffices to establish

$$(4.1) \qquad\qquad \log \mathbb{P}(G_1) \geq -(3 + o(1))\log K.$$

Let $a_l := F^{(K+1)}(\sum_{i=l}^{l+K} Y_i)$ for each $l \in [1, K+1]$. Then

$$\mathbb{P}(G_1) = \mathbb{E}\mathbb{P}(G_1 | \mathcal{F}) = \mathbb{E}\prod_{l=1}^{K+1} a_l.$$

If $a_l \geq 1 - 1/K$ for each $l \in [1, K+1]$, then $\prod_{l=0}^{K} a_l \geq e^{-1} + o(1)$, so (4.1) follows from

$$(4.2) \qquad\qquad \mathbb{P}\left(\min\{a_l : 1 \leq l \leq K+1\} \geq 1 - \frac{1}{K}\right) \geq cK^{-3}.$$

Let $\mathcal{F}^*$ be the $\sigma$-field generated by the unordered pair of sets $\{Y_1, \ldots, Y_{K+1}\}$ and $\{Y_{K+2}, \ldots, Y_{2K+2}\}$. Then $\min\{a_1, a_{K+2}\} \in \mathcal{F}^*$. Furthermore, conditional on $\mathcal{F}^*$, the collection $\{S_l := \sum_{i=1}^{l}(Y_{i+K+1} - Y_i) : 1 \leq l \leq K+1\}$ has exchangeable increments (generated by continuous distribution i.i.d. picks, so ties in the partial sum sequence $S_\cdot$ happen with probability 0) that are symmetric about 0. Now note the following consequence of exchangeability: Conditioned on all the increments, if their total sum is positive, then the probability that the minimum occurs at the beginning, that is, all the intermediate sums are positive, is at least $1/K$. Namely, all cyclic permutations of the increments are equally distributed and almost surely there is at least one such permutation for which the minimum is achieved at step 0.

Therefore,

$$\mathbb{P}(\min\{S_l : 1 \leq l \leq K+1\} > 0) \geq \tfrac{1}{2}K^{-1}.$$

When $\{\min\{S_l : 1 \leq l \leq K+1\} > 0\}$ occurs, we have $\min\{a_l : 1 \leq l \leq K+1\} = a_1$. Hence, by conditioning on $\mathcal{F}^*$ first, the probability on the left-hand side of (4.2) is at least

$$\tfrac{1}{2}K^{-1}\mathbb{P}(\min\{a_1, a_{K+2}\} \geq 1 - 1/K),$$

and by independence of $a_1$ and $a_{K+2}$ (recall that $N > 2K + 1$) this is equal to

$$(\tfrac{1}{2} + o(1))K^{-1}K^{-2},$$

proving (2.4). □

The proofs of Theorems 2.3 and 2.4 are similar to the argument used to prove the second inequality of Theorem 2.2. Having specific distributions to work with makes the arguments simpler and the results sharper (cf. Conjecture 1).

PROOF OF THEOREM 2.3. Cover the interval $[N] := \{1, \ldots, N\}$ with $L := \lceil N/((1-y)(1+K)) \rceil$ intervals of size $\lceil (1-y)(K+1) \rceil$. Denote these intervals by $I_1, \ldots, I_L$. Positive association again implies that

$$\mathbb{P}(\mathcal{H}') \geq \mathbb{P}(\mathcal{H}'_j \ \forall j \in I_1)^L.$$

Let $I'$ denote the interval of length $\lfloor y(K+1) \rfloor$ adjacent to and just preceding $I_1$. If the maximum of the collection $\{Y_j, Y_{l,i} : j \in I', l \in I_1, 0 \leq i \leq K\}$ is $Y_{j_0}$ for some $j_0 \in I'$, it follows that $\mathcal{H}'_j$ occurs for each $j \in I_1$. The last claim follows directly from definition (2.5) since for such $j_0$ we have $Y_{j_0} \leq \max_{i=j-K}^{j} Y_j$ whenever $j \in I_1$.

The probability of

$$\Big\{\max_{j \in I'} Y_j > \max_{l \in I_1, 0 \leq i \leq K} Y_{l,i}\Big\},$$

up to corrections for integer roundoff, is clearly equal to $y(K+1)/[(1-y)(K+1)^2 + y(K+1)]$. Thus

$$\mathbb{P}(\mathcal{H}') \geq \Big[(1+o(1))\frac{y(K+1)}{(1-y)(K+1)^2 + y(K+1)}\Big]^L.$$

Choosing $y = y(K) = 1/\log K$ as before suffices to prove the theorem. $\square$

PROOF OF THEOREM 2.4.   Keeping the notation from the previous proof, we need to estimate $\mathbb{P}(\mathcal{H}_j \forall j \in I_1)$ when $F$ is the Cauchy distribution. Define events:

(i) $A := \{\max_{j \in I'} Y_j \geq 2(K+1)^2\}$;

(ii) $B := \{\sum_{j \in I_1 \cup I_0} 0 \vee (-Y_j) < (K+1)^2\}$;

(iii) $C := \{\max_{j \in I_1} \sum_{i=0}^K Y_{j,i} \leq (K+1)^2\}$.

Here $I_0$ is the interval of length $K$ preceding $I_1$ so that for $y(K) < 1$ (which will be the case) $I' \subset I_0$. Note that $A \cap B \cap C \subset \bigcap_{j \in I_1} \mathcal{H}_j$ since on $A \cap B \cap C$ we have, if $j \in I_1$, both

$$\left\{ \sum_{i=j}^{j-K} Y_i \geq \max_{j \in I'} Y_j + \sum_{j \in I_1 \cup I_0} (0 \vee (-Y_j)) > (K+1)^2 \right\}$$

and

$$\left\{ \sum_{i=0}^K Y_{j,i} \leq (K+1)^2 \right\}.$$

It is not difficult to check that

$$P(A) = ((2\pi)^{-1} + o(1)) y K^{-1},$$

$$P(B) \geq \exp\left( \frac{-(2-y)^2}{\pi} \right) + o(1),$$

$$P(C) = \exp\left( \frac{-(1-y)}{\pi} \right) + o(1).$$

Indeed,

$$P(A) = 1 - P\left( \max_{j \in I'} Y_j \leq 2(K+1)^2 \right)$$

$$= 1 - \left( 1 - \frac{(1/\pi + o(1))}{2(K+1)^2} \right)^{y(K+1)} = -\frac{y}{K+1} \left( \frac{1}{2\pi} + o(1) \right),$$

$$P(B) \geq P\left( \max_{j \in I_1 \cup I_0} 0 \vee (-Y_j) < \frac{K+1}{2-y} \right)$$

$$= \left( 1 - \frac{2-y}{\pi(K+1)} \right)^{(K+1)(2-y)} \geq \exp\left( -\frac{(2-y)^2}{\pi} \right) + o(1)$$

and

$$P(C) = P\left( \sum_{i=0}^K Y_{1,i} \leq (K+1)^2 \right)^{(1-y)(K+1)} = \left( \int_{-\infty}^{K+1} \frac{1}{\pi(1+y^2)} \, dy \right)^{(1-y)(K+1)}$$

$$= \left(1 - \frac{(1/\pi + o(1))}{K+1}\right)^{(1-y)(K+1)} = \exp\left(-\frac{1}{\pi}(1-y)\right) + o(1).$$

Another application of positive association shows that

$$P(A \cap B \cap C) \geq \left(\frac{e^{-5/\pi}}{2\pi} + o(1)\right)\frac{y}{K}$$

so that

$$\mathbb{P}(\mathcal{H}) \geq \left[\left(\frac{e^{-5/\pi}}{2\pi} + o(1)\right)\frac{y}{K}\right]^{L},$$

and taking the logarithm, with $y(K) = \log(K)^{-1}$, completes the proof. $\quad\square$

## 5. The fat tail when $N/K$ remains bounded.
This section provides a proof of Theorem 2.5. In particular, in this section we derive asymptotic formulae for $p_{\mathrm{fat}}(N,K)$ that are valid as $N, K \to \infty$, uniformly as long as $N/K$ remains bounded. Probability estimates come from the following algorithm for checking whether $\mathcal{H}'$ has occurred.

1. Initialize $r = 1$ and $\mathcal{C}$ to be the collection of variables $\{Y_j, Y_{j,i} : 1 \leq j \leq N, 0 \leq i \leq K\}$.
2. Find the maximum of the variables in $\mathcal{C}$.
3.   (a) If this maximum is one of the variables $Y_{j,i}$, then output FALSE and stop.

    (b) Else, let $j_r$ be the index such that the maximum occurred at $Y_{j_r}$.
4. Remove from $\mathcal{C}$ the variables $Y_{j,i}$ for $j_1 \leq j \leq j_1 + K$, $0 \leq i \leq K$ (these are no longer relevant since no matter what their value is, we know that $\bigcap_{j=j_1}^{j_1+K} \mathcal{H}'_j$ has occurred, and other $\mathcal{H}'_l$'s do not depend on the values of $Y_{j,i}$, $j_1 \leq j \leq j_1 + K$, $0 \leq i \leq K$, anyhow), and also remove the variable $Y_{j_r}$.
5.   (a) If the collection $\mathcal{C}$ contains no more variables $Y_{j,i}$, then output TRUE and stop.

    (b) Else, set $r$ to $r + 1$ and go to Step 2.

Clearly $\mathcal{H}' = \{\text{algorithm stops at TRUE}\}$. We may think of the output as containing all values of $j_r$ found before stopping, so that in addition to the indicator function of the event $\mathcal{H}'$, the algorithm outputs the random variables $R, j_1, \ldots, j_R$, where $R$ is the maximum value for which the first Step 3(b) (the else statement) is executed. Recall that $r_0 := \lceil N/(K+1)\rceil$ is a lower bound for $R$, provided the output is TRUE. The possible values for the sequence $\mathbf{j}$ when it is of length $R = r$ are precisely the set $S(r)$ of sequences that satisfy both of the following statements:

($*$) For every $i \in [N]$ there is an $s \leq r$ for which $0 \leq i - j_s \leq K$.
($**$) No initial segment of $\mathbf{j}$ satisfies property ($*$).

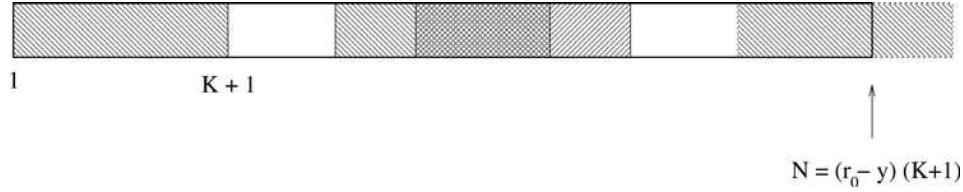

Letting $\mathcal{H}(\mathbf{j})$ denote the event that $\mathcal{H}'$ occurs and the algorithm outputs the witnessing sequence $\mathbf{j}$, we may decompose $\mathcal{H}'$ into a disjoint union by setting $\mathcal{H}(r) := \bigcup_{\mathbf{j} \in S(r)} \mathcal{H}(\mathbf{j})$ and

$$\mathcal{H}' = \bigcup_r \mathcal{H}(r) = \bigcup_r \bigcup_{\mathbf{j} \in S(r)} \mathcal{H}(\mathbf{j}).$$

Given $1 \le s \le r + 1$ and any sequence $\mathbf{j}$ of length $r$ containing distinct elements of $[N]$, define

$$\mathsf{missed}(s, \mathbf{j}) := \{j \in [N] : j - j_t \notin \{0, \ldots, K\} \ \forall \, t < s\},$$

$$M(s, \mathbf{j}) := |\mathsf{missed}(s, \mathbf{j})|.$$

Vacuously, $M(1, \mathbf{j}) = N$ for all $\mathbf{j}$. Figure 1 illustrates this definition when $r_0 = 4$. In the illustration, the intervals $[j_s, \ldots, j_s + K]$ are shaded, $j_1$ is equal to $K + 1$, one interval overlaps with $[1, K + 1]$ modulo $N$ and the other two intervals also overlap. Figure 1 also illustrates a general fact, namely that the set $\mathsf{missed}(s, \mathbf{j})$ (the white space between the shaded intervals) is always composed of no more than $s$ intervals (i.e., the unshaded set has at most $s$ connected pieces), where adjacent white intervals are separated by a distance of at least $K + 1$.

One further observation is that for all $s$ and $\mathbf{j}$,

$$(5.1) \qquad N \ge M(s, \mathbf{j}) \ge N - (s - 1)(K + 1).$$

Conditional on the event $R \ge r + 1$ and on $j_1, \ldots, j_r$, the values of the variables remaining in $\mathcal{C}$ at stage $r$ are i.i.d., so the conditional probability of $j_{r+1} = j$ for any $j \notin \{j_1, \ldots, j_r\}$ is equal to the reciprocal of the number of variables remaining in $\mathcal{C}$, that is, $1/(N - r + (K + 1)M(r, \mathbf{j}))$. Applying this inductively yields

$$(5.2) \qquad \mathbb{P}(\mathcal{H}(\mathbf{j})) = \prod_{s=1}^{r} \frac{1}{N - (s - 1) + (K + 1)M(s, \mathbf{j})}.$$

The $(K + 1)M(s, \mathbf{j})$ contribution above comes from the number of $Y_{l,i}$ variables that are still in $\mathcal{C}$. The above computation can be generalized in the following useful way. Define the event $\mathcal{H}^*(\mathbf{j})$ by

$$\mathcal{H}^*(\mathbf{j}) := \mathcal{H}' \cap \{\mathbf{j} \text{ is an initial segment of the output of the algorithm}\}.$$

When $\mathbf{j}$ of length $r$ is an element of $S(r)$, $\mathcal{H}(\mathbf{j}) = \mathcal{H}^*(\mathbf{j})$; otherwise $\mathcal{H}(\mathbf{j})$ is empty and the right-hand side in (5.2) computes the probability of outputting $\mathbf{j}$ as an initial segment. To obtain $\mathbb{P}(\mathcal{H}^*(\mathbf{j}))$ from this, one must multiply the right-hand side in (5.2) by the probability $Q(\mathbf{j})$ that, conditional on the initial segment being $\mathbf{j}$, the algorithm eventually outputs TRUE. We compute an upper bound on $Q(\mathbf{j})$, for $\mathbf{j}$ of length $r$, as follows. For each interval $I = [a, b] \subseteq \mathsf{missed}(r, \mathbf{j})$, for the $\mathcal{H}^*(\mathbf{j})$ to happen, it is necessary that $\max_{a-K \leq j \leq b} Y_j$ be greater than $\max_{j \in I, 0 \leq i \leq K} Y_{j,i}$. This probability of $\{\max_{a-K \leq j \leq b} Y_j > \max_{j \in I, 0 \leq i \leq K} Y_{j,i}\}$ equals

$$\frac{b+K+1-a}{b+K+1-a+(b+1-a)K} = \frac{(b+1-a)+K}{b+1-a+(b+2-a)K}$$

$$\leq \frac{1}{K+1} + \frac{1}{b+2-a}.$$

If $\mathsf{missed}(r, \mathbf{j})$ is composed of more than one interval, the probabilities for each interval are multiplied (since they are at least $K+1$ units apart, everything is independent) and, therefore, for a given $M(r, \mathbf{j})$, the upper bound on $Q(\mathbf{j})$ is greatest when $\mathsf{missed}$ has only one interval and we may take as an upper bound

$$(5.3) \qquad Q(\mathbf{j}) \leq \frac{1}{M(r, \mathbf{j})} + \frac{1}{K}.$$

We now bound the number of sequences $\mathbf{j}$ that produce a given value of $M(r, \mathbf{j})$.

LEMMA 5.1. *Let $N = (r - y)(K + 1)$. Then the number of sequences $\mathbf{j}$ of length $r$ with $M(r+1, \mathbf{j}) = \iota$ is at most*

$$NC(r)(yK + \iota)^{r-2}.$$

PROOF. By symmetry, it suffices to consider only sequences for which $j_1 < \cdots < j_r$ in cyclic order modulo $N$ and then multiply by $(r-1)!$. By convention, we let $j_0 := j_r - N$. For $1 \leq s \leq r$, consider the quantities $A_s := j_{s-1} + K + 1 - j_s$ to be unknown and satisfying the following two nice properties:

(a) $\sum_{s=1}^{r} A_s = j_0 - j_r + r(K+1) = -N + r(K+1) = y(K+1)$;
(b) $\sum_{s=1}^{r} (-A_s) \vee 0 = \sum_{s=1}^{r} [(j_s - j_{s-1}) - (K+1)] \wedge 0 = \iota$.

Property (b) is a consequence of the fact that the length of the unique (white) interval that contributes to $\mathsf{missed}(r+1, \mathbf{j})$, which is contained in $[j_{s-1}, j_s]$, equals $[(j_s - j_{s-1}) - (K+1)] \wedge 0$. The sequence $(A_1, \ldots, A_r)$ and the value $j_1$ together determine $\mathbf{j}$. The number of possible sequences $(A_1, \ldots, A_r)$ above may be bounded as follows. Let $S_+$ be the set of indices $i$ for which $A_i \geq 0$. Given $S_+$, the subsequence $(A_i : i \in S_+)$ is a sequence of nonnegative integers that sum to $y(K+1) + \iota$. These sequences are called *compositions* of $y(K+1) + \iota$ into $|S_+|$ parts, and the number of such compositions is $\binom{y(K+1)+\iota+|S_+|-1}{|S_+|-1}$ ([10], page 14). Similarly, $(A_i : i \notin S^+)$ is a composition of $\iota$ into $r - |S_+|$ parts, and the number of these is $\binom{\iota+r-|S_+|-1}{r-|S_+|-1}$. We claim that the product of the above two binomial coefficients is bounded above by $C_0(r)(y(K+1) + \iota)^{r-2}$. Indeed, the product equals

$$\frac{(y(K+1)+\iota+|S_+|-1)!}{(|S_+|-1)!(y(K+1)+\iota)!} \cdot \frac{(\iota+r-|S_+|-1)!}{(\iota)!(r-|S_+|-1)!}.$$

Clearly $|S_+| \leq y(K+1) + \iota$ and $r - |S_+| \leq \iota$, which implies

$$\frac{(y(K+1)+\iota+|S_+|-1)!}{(y(K+1)+\iota)!} \leq [2(y(K+1)+\iota)]^{|S_+|-1}$$

and

$$\frac{(\iota+r-|S_+|-1)!}{(\iota)!} \leq [2(y(K+1)+\iota)]^{r-|S_+|-1}.$$

Thus, for a given $S_+$, there are at most $NC_0(r)(y(K+1) + \iota)^{r-2}$ such $\mathbf{j}$ sequences ($N$ comes from the choice of $j_1$). Summing over at most $2^r - 2$ values of $S_+$ proves the lemma. $\quad\square$

As mentioned prior to the statement of Theorem 2.6, the complexity in the behavior of $p_{\text{fat}}(N, K)$ is due to transitions in the number of $Y_j$ variables with large values from one integer to the next higher. We separate the argument into several cases, the first three being restricted to $r_0 = \lceil N/K \rceil \geq 3$:

1. $r_0 - 1 + \varepsilon \leq N/(K+1) \leq r_0 - \varepsilon$;
2. $r_0 - \varepsilon \leq N/(K+1) \leq r_0$;
3. $r_0 - 1 \leq N/(K+1) \leq r_0 - 1 + \varepsilon$;
4. $r_0 = 2$.

The analyses of Cases 2 and 3 actually cover Case 1 since one could take $\varepsilon = 1/2$, but since the argument is easier for values of $N/(K+1)$ not too close to an integer, we prefer to present this as the first case.

CASE 1.   We first compute $\mathbb{P}(\mathcal{H}(r_0))$. For each $\mathbf{j} \in S(r_0)$ and each $s \leq r_0$, the expression (5.2) and the bounds (5.1) imply

$$c(\varepsilon) \leq K^{2r_0}\mathbb{P}(\mathcal{H}(\mathbf{j})) \leq C(\varepsilon) \qquad (\text{recall } N \sim r_0 K).$$

Together with the fact that $S(r_0)$ has cardinality $\Theta(K^{r_0})$ (see below for details), this immediately implies that

$$\mathbb{P}(\mathcal{H}(r_0)) = \Theta(K^{-r_0}).$$

In this case, we claim that $\mathbb{P}(\mathcal{H}(r))$ is maximized at $r = r_0$. With $\mathbb{P}(\mathcal{H}(r-1))$ trivially being zero, this statement and the theorem follow from a more precise estimate of $\mathbb{P}(\mathcal{H}(r_0))$ and a bound on $\mathbb{P}(\mathcal{H}^*(r_0+1))$.

Let $\mathsf{T}$ be the $r_0$-dimensional torus of $r_0$-tuples in $\mathbb{R}/\mathbb{Z}$, with addition modulo 1 and unit Lebesgue measure $\lambda$. For $y \in [0, r_0]$, define a subset $\mathsf{T}(y) = \mathsf{T}(y, r_0) \subseteq \mathsf{T}$ to be the set of $\mathbf{x} = (x_1, \ldots, x_{r_0})$ such that for all $z$ there is a $j \le r_0$ with $x_j - 1/(r_0 - y) \le z \le x_j$. Consider the mapping of $S(r_0)$ into $\mathsf{T}$ by

$$(5.4) \qquad\qquad \mathbf{j} \mapsto \mathbf{j}/N.$$

The set $S(r_0)$ then maps into the set $\mathsf{T}(y)$ for $y = r_0 - (N/(K+1))$. In fact, for any $U \subseteq \mathsf{T}(y)$, the cardinality of the subset of $S(r_0)$ that maps into $U$ under (5.4) is equal to $(1 + o(1))N^{r_0}\lambda(U)$ uniformly in $N/K$ as $N \to \infty$. Furthermore, for $\mathbf{j} \in S(r_0)$,

$$(5.5) \qquad\qquad \mathbb{P}(\mathcal{H}(\mathbf{j})) = (1 + o(1))N^{-r_0}K^{-r_0}\eta\left(\frac{\mathbf{j}}{N}\right),$$

where

$$(5.6) \qquad\qquad \eta(\mathbf{x}) = \prod_{s=1}^{r_0} \frac{1}{\tilde{M}(s, \mathbf{x})}$$

and $\tilde{M}(s, \mathbf{x})$ is the measure of $[0, 1] \setminus \bigcup_{t=1}^s [x_t - K/N, x_t]$. Let $y = r_0 - N/(K+1)$ and note that $y$ equals $j/(K+1)$ when $N = r_0(K+1) - j$ for $j > 0$. By bounded convergence, we then have

$$(5.7) \qquad\qquad K^{r_0}\mathbb{P}(\mathcal{H}(r_0)) \to f_{r_0}(y) := \int_{\mathsf{T}(y \wedge 1; r_0)} \eta(\mathbf{x}) \, d\lambda(\mathbf{x})$$

(note here that since $y \in [\varepsilon, 1 - \varepsilon]$, $y \wedge 1 = y$) as $N \to \infty$, uniformly in $N/K$, with $f_{r_0}(\cdot)$ bounded, continuous and nondecreasing. This is the $f_r$ term in the last line of Table 1.

Next, we compute an upper bound for the event $\mathcal{H}^*(r_0+1) := \bigcup\{\mathcal{H}^*(\mathbf{j}) : \mathbf{j} \notin S(r_0), |\mathbf{j}| = r_0\}$ that an output of TRUE requires at least $r_0 + 1$ covering intervals. Multiplying the right-hand side of (5.2) by $Q(r_0, \mathbf{j}) = Q(\mathbf{j})$, using (5.3) with $r = r_0$ and using the fact that $M(s, \mathbf{j}) \ge C(\varepsilon)K$ for $s \le r_0$, we see that

$$\mathbb{P}(\mathcal{H}^*(r_0+1)) \le \sum_{\mathbf{j} \notin S(r_0)} Q(\mathbf{j}) \prod_{s=1}^{r_0} \frac{1}{N - (s-1) + (K+1)M(s, \mathbf{j})}$$

$$\le \sum_{s=1}^N \sum_{M(r_0+1, \mathbf{j}) = s} C\frac{1}{s}\frac{1}{K^{2r_0}},$$

where $C$ represents a constant that depends only on $r_0$ and $\varepsilon$, and the sum is over sequences $\mathbf{j}$ of length $r_0$. By Lemma 5.1, we may further bound this from above by

$$\mathbb{P}(\mathcal{H}^*(r_0+1)) \leq \sum_{s=1}^{N} C(yK+s)^{r_0-2}\frac{1}{s}\frac{N}{K^{2r_0}}\left(\frac{N}{K} \leq r_0+1\right)$$

$$(5.8) \qquad \leq C'\frac{1}{K^{r_0+1}}\sum_{s=1}^{N}\frac{1}{s}\left(y+\frac{s}{K}\right)^{r_0-1}\left(\sum_{s=1}^{N}\frac{1}{s} \sim \log(N) \sim \log(K)+c\right)$$

$$\leq C''\frac{\log K}{K^{r_0+1}}.$$

Together with (5.7), this establishes that

$$(5.9) \qquad\qquad K^{r_0}\mathbb{P}(\mathcal{H}') \to f_{r_0}(y).$$

When $\varepsilon < y < 1-\varepsilon$, the term containing $f_r$ in the last line of Table 1 dominates the term containing $f_{r+1}$ since $f_r(\varepsilon) > 0$, so this proves the theorem in the case $\varepsilon < y < 1-\varepsilon$ and $N/(K+1) \geq 3$.

Case 2. This is quite similar to the previous case. The part where we estimated (5.7) goes through unchanged, only now $f_{r_0}$ tends to zero as $N/(K+1) \to r_0^-$ and we need to find the asymptotic rate to compare to the $f_{r_0+1}$ term.

Lemma 5.2. *The measure $\lambda(\mathsf{T}(y,r_0))$ of $\mathsf{T}(y,r_0)$ is asymptotically $y^{r_0-1}/(r_0-y)^{r_0-1}$ near $y=0$.*

Proof. The set $\mathsf{T}(y)$ is invariant under translation of each coordinate by a constant, so by symmetry the measure is the same as the $(r_0-1)$-dimensional measure of the fiber of $\mathsf{T}(y)$, where $x_1 = 0$. By permutation invariance, this is equal to $(r_0-1)!$ times the measure of the subset of $\mathsf{T}(y)$, where $0 = x_1 < x_2 \cdots < x_{r_0}$. Such a point is in $\mathsf{T}(y)$ if and only if the quantities $x_i + 1/(r_0-y) - x_{i+1}$, for $1 \leq i \leq r_0-1$, are positive numbers summing to at most $y/(r_0-y)$. In fact, the mapping that maps each $\mathbf{x}$ in the fiber to the sequence $(x_1 + K/N - x_2, \ldots, x_{r_0-1} + K/N - x_{r_0})$ is an isometry. The $(r_0-1)$-dimensional simplex of positive numbers summing to at most $y/(r_0-y)$ has volume $y^{r_0-1}/((r_0-y)^{r_0-1}(r_0-1)!)$, which proves the lemma. □

As $y \to 0$, the factors $1/\tilde{M}(s,\mathbf{x})$ converge to $r_0/(r_0-(s-1))$, since the only way for a vector to be in $\mathsf{T}(y)$ is for it to have $r_0$ approximately evenly spaced coordinates. Therefore, the function $\eta$ defined in (5.6) converges to the constant $r_0^{r_0}/r_0!$ on $\mathsf{T}(y)$, and we have

$$f_{r_0}(y) = \int_{\mathsf{T}(y)} \eta(\mathbf{x})\,d\lambda(\mathbf{x}) \sim \frac{r_0^{r_0}}{r_0!}\lambda(\mathsf{T}(y)) \sim \frac{y^{r_0-1}}{r_0!}.$$

Since the contribution of $\mathbb{P}(\mathcal{H}(r_0+1))$ to $\mathbb{P}(\mathcal{H}')$ is no longer negligible, we must compute it a little more precisely as well. If we write it as an integral analogous to (5.7), we find, for $r_0 \geq 3$, that the integral $\int_{\mathsf{T}(1)} \eta(\mathbf{x})\,d\lambda(x)$ exists as an improper integral, but the integral over $\mathsf{T}(y)$ diverges for $y > 1$. We have shown that $K^{r_0}\mathbb{P}(\mathcal{H}(r_0)) \sim y^{r_0-1}/(r_0-1)!$ as $y \to 0$, and we have an upper bound (5.8) on $\mathbb{P}(\mathcal{H}^*(r_0+1))$. When $y \geq K^{-1/r_0}$, these two together show that still

$$K^{r_0}\mathbb{P}(\mathcal{H}') \sim f_{r_0}(y).$$

Assume therefore that

$$(5.10) \qquad\qquad\qquad y \leq K^{-1/r_0}.$$

We cannot immediately conclude for $0 \leq y \leq K^{-1/r_0}$ that

$$K^{r_0+1}\mathbb{P}(\mathcal{H}^*(r_0+1)) \to f_{r_0+1}(1)$$

and it is our remaining task to verify the above statement. One part of this is easy. For any positive $L$, the function $\eta\mathbf{1}_{\eta>L}$ is bounded and, as $L \to \infty$, these functions converge in $L^1$ to $\eta$ as long as $\eta \in L^1$, which is the case since we have assumed that $r_0 \geq 3$. Equivalently, the function

$$g(L) := \int_{\mathsf{T}(1;r_0+1)} \eta(\mathbf{x})\mathbf{1}_{\eta(\mathbf{x}) \geq L}\,d\lambda(\mathbf{x})$$

converges to 0 as $L \to \infty$ and, by bounded convergence, we may approximate the truncated sum of the terms in (5.5) by a truncated integral as $K \to \infty$:

$$(5.11) \quad K^{r_0+1}\mathbb{P}(\mathcal{H}(r_0+1) \cap \{\eta(\mathbf{j}/N) \leq L\}) \to (1-g(L))f_{r_0+1}(1).$$

The theorem, in Case 2, follows if we can show that

$$(5.12) \qquad \mathbb{P}\Big(\eta\Big(\frac{\mathbf{j}}{N}\Big) \leq L, \mathcal{H}^*(r_0+2)\Big) \leq C(L)\frac{1}{K^{r_0+2}},$$

$$(5.13) \qquad \mathbb{P}\Big(\eta\Big(\frac{\mathbf{j}}{N}\Big) \geq L, \mathcal{H}^*(r_0+1)\Big) \leq c(L)\frac{1}{K^{r_0+1}}$$

for $c(L) \to 0$ as $L \to \infty$, uniformly in $K$. Indeed if these two hold, then for $L$ large enough so that $c(L) < \delta/2$ and $K$ then chosen large enough so that $C(L)/K < \delta/2$, we have

$$\mathbb{P}\Big[\mathcal{H}' \Big\backslash \Big(\mathcal{H}(r_0+1) \cap \Big\{\eta\Big(\frac{\mathbf{j}}{N}\Big) \leq L\Big\}\Big)\Big] \leq \frac{\delta}{K^{r_0+1}},$$

which together with (5.11) finishes Case 2.

To prove (5.12), we may use the same argument that proved (5.8), but with $r_0$ replaced by $r_0 + 1$. We sum over sequences $\mathbf{j}$ of length $r_0 + 1$ to get

$$\mathbb{P}\left(\eta\left(\frac{\mathbf{j}}{N}\right) \leq L, \mathcal{H}^*(r_0 + 2)\right) \leq \sum_{\mathbf{j} \notin S(r_0+1)} Q(\mathbf{j}) \prod_{s=1}^{r_0+1} \frac{1}{N - (s-1) + (K+1)M(s,\mathbf{j})}$$

$$\leq \sum_{s=1}^{N} \sum_{M(r_0+1,\mathbf{j})=s} C(L)\frac{1}{s}\frac{1}{K^{2r_0+2}}.$$

Here we have used the fact that $\eta(\mathbf{j}/N) \leq L$ to bound the product in the first line by $C(L)K^{-2r_0-2}$; equation (5.3) is valid for any $r$, so there is no trouble replacing $r_0$ by $r_0 + 1$ here. At the next step, instead of requiring Lemma 5.1, we require only the trivial bound on the number of sequences $\mathbf{j}$ of length $r_0 + 1$ with $M(r_0 + 2, \mathbf{j}) = j$, namely $CK^{r_0}$. Following the path to (5.8) leads this time to (5.12).

To prove (5.13), observe first that $\eta(\mathbf{j}/N) \geq L$ implies $M(r_0 + 1, \mathbf{j}) \leq \varepsilon(L)K$ for some function $\varepsilon(L)$ going to zero as $L \to \infty$. This follows from expression (5.2), according to which all the factors $1/\tilde{M}(s,\mathbf{j})$ in the definition of $\eta$ are bounded from below except for the factor with $s = r_0 + 1$, which is of order $K/M(r_0 + 1, \mathbf{j})$. Hence,

$$\mathbb{P}\left(\eta\left(\frac{\mathbf{j}}{N}\right) \geq L, \mathcal{H}^*(r_0 + 1)\right)$$

$$\leq \sum_{s=1}^{\varepsilon(L)K} \sum_{\mathbf{j}:M(r_0+1,\mathbf{j})=s} Q(r_0,\mathbf{j}) \prod_{t=1}^{r_0} \frac{1}{N - (t-1) + (K+1)M(t,\mathbf{j})}$$

$$\leq \sum_{s=1}^{\varepsilon(L)K} C(r_0,\varepsilon)(N - r_0 K + 2s)^{r_0-1} \frac{1}{s} K^{-2r_0}$$

$$\leq \frac{C(r_0,\varepsilon)}{K^{r_0+1}} \sum_{s=1}^{\varepsilon(L)K} \left(\frac{N}{K} - r_0 + 2\frac{s}{K}\right)^{r_0-1}\frac{1}{s}.$$

This sum is at most twice the integral for which it is an upper Riemann. To be precise, we consider the sum as a step function, change variables to $x = (s+1)/K$, and compare the upper and lower Riemann sums to integrals, concluding that

$$K^{r_0+1}\mathbb{P}\left(\eta\left(\frac{\mathbf{j}}{N}\right) \geq L, \mathcal{H}^*(r_0 + 1)\right)$$

$$\leq 2C(r_0,\varepsilon)\int_0^{\varepsilon(L)}\left(\frac{N}{K} - r_0 + 2x\right)^{r_0}\frac{1}{x + K^{-1}}\,dx.$$

As a family of functions on $[0, 1]$, the integrands form a uniformly integrable family as long as $N/K - r_0 \leq K^{-\alpha}$ for some $K$. By assumption (5.10), this inequality is indeed satisfied, and we may conclude that the integral from 0 to $\varepsilon(L)$ tends to zero uniformly in $K$ as $\varepsilon(L) \to 0$. This finishes the proof of (5.13) and therefore of Case 2. We go onto Case 4, coming back to Case 3 later since it uses some of the computations from Case 4.

CASE 4. When $r_0 = 2$, the computation is particularly simple without using the continuous approximation. The first term in the product in (5.2) is always $1/(N(K+2))$. By symmetry,

$$\mathbb{P}(\mathcal{H}(2)) = N \sum_{\substack{\mathbf{j} \in S(2) \\ j_1 = 1}} \mathbb{P}(\mathcal{H}(\mathbf{j})).$$

For $j_1 = 1$, so that $\mathbf{j}$ satisfies property $(*)$, it is necessary to choose $N - K \leq j_2 \leq K + 2$. Also, if $j_1 = 1$, then $\mathsf{missed}(\mathbf{j}, 2) = \{K + 2, \ldots, N\}$ and the second factor in (5.2) is always $1/(N - 1 + (K + 1)(N - K - 1))$. Thus, letting $j = 2(K + 1) - N \in \{0, \ldots, K\}$, we have

$$
\begin{aligned}
(5.14) \quad p(\mathcal{H}(2)) &= \frac{N}{N(K+2)} \frac{2K + 3 - N}{N - 1 + (K + 1)(N - K - 1)} \\
&= \frac{1}{K+2} \frac{j + 1}{N - 1 + (N - K - 1)(K + 1)} \\
&= (1 + o(1)) \frac{1}{K^2} \frac{j + 1}{N - (K + 1) + ((N - 1)/(K + 1))} \\
&= (1 + o(1)) \frac{1}{K^2} \frac{j + 1}{K + 3 - j}.
\end{aligned}
$$

For $\mathbb{P}(\mathcal{H}(3))$, a similarly direct argument ensues. If $\mathcal{H}(3)$ occurs via $\mathcal{H}(\mathbf{j})$ for some $\mathbf{j} \in S(3)$ with $j_1 = 1$, then since $\mathcal{H}(2)$ does not occur, either $j_2 \in [2, N - (K + 1)]$ or $j_2 \in [K + 3, N]$. In the former case, $j_3 \in [N - K, j_2 + K + 1]$, while in the latter case, $j_3 \in [j_2 - K - 1, K + 2]$. For the first of the two cases, we then have a contribution to $p(\mathcal{H}(3))$ of

$$
\begin{aligned}
\frac{1}{K+2} &\sum_{j_2 = 2}^{N - K - 1} \sum_{j_3 = N - K}^{j_2 + K + 1} \frac{1}{N - 1 + (K + 1 - j)(K + 1)} \\
&\qquad\qquad \times \frac{1}{(N - 2) + (N - j_2 - K)(K + 1)} \\
(5.15) \quad &= \frac{1}{K + 2} \frac{1}{N - 1 + (K + 1 - j)(K + 1)}
\end{aligned}
$$

$$\times \sum_{j_2=2}^{N-K-1} \frac{j+j_2}{N-2+(K+1)(N-K-j_2)}$$

$$= (1+o(1)) \frac{1}{K^3(1-j/K)} \sum_{s=1}^{N-K-2} \frac{1-s/K}{s}$$

$$= (1+o(1)) \frac{\log K}{K^3(1-j/K)}.$$

Here the third equality comes from the substitution $s = N - K - j_2$ and the definition of $j$ as $2(K+1) - N$, while the $(1+o(1))$ term comes from factors of order $(1 + O(1/K))$ that remain once we remove three factors of $K$ from the top and bottom of the fraction preceding the sum and one factor of $K$ from the top and bottom of the summand. The computation for the second case is symmetrical, leading to

$$(5.16) \qquad p(\mathcal{H}(3)) = (2+o(1)) \frac{\log K}{K^3(1-j/K)}.$$

Comparing (5.16) to (5.14), we see that the former is dominant when $j = o(\log K)$, the latter when $\log K = o(j)$ and both contribute when $j = \Theta(\log K)$. In particular, (5.14) contributes only when $j \to \infty$, in which case the contribution is $(1+o(1))\frac{1}{K}(\frac{1}{K+3-j} - \frac{1}{K})$, while (5.16) contributes only when $j = o(K)$, in which case the contribution is $(2+o(1))\log K/K^3$. From these the first line in the Table 1 follows as a lower bound, with an identical upper bound yet to follow if we show that changing $\mathcal{H}(3)$ to $\mathcal{H}^*(4)$ produces no change to the asymptotics.

The difference between $\mathcal{H}(3)$ and $\mathcal{H}^*(4)$ is that in the latter case, $j_3$ can be element of $\mathsf{missed}(3, (j_1, j_2))$. These are all $j'$ not in the interval $[1, j_2]$, so the numerator $j + j_2$ of (5.15) becomes $N - j_2$. This changes the $1 - s/K$ in the numerator of the subsequent line to $1 + s/K$, which does not affect the sum asymptotically since all the contribution come from $s = o(K)$.

CASE 3. The analysis of the $\mathbb{P}(\mathcal{H}(r_0 + 1))$ term in Case 2 works just as well for $N$ slightly greater than $r_0(K+1)$, and this becomes the $f_r$ term in the last line of the table for $r = r_0 + 1$. Since $r_0 \geq 3$, Case 2 handles the $f_r$ terms for $r \geq 4$. It remains only to analyze the $f_3$ term appearing in line 2 of the table.

We borrow the analysis from Case 4. Now the event $\mathcal{H}(2)$ cannot happen, so we need to evaluate $\mathbb{P}(\mathcal{H}(3))$, show it gives the asymptotics stated in the theorem and then show that adding $\mathbb{P}(\mathcal{H}^*(4))$ does not alter the asymptotics. Let $N = 2(K+1) + j$. Assume $j_1 = 1$, so the first interval thrown out of $\mathcal{C}$ is $[1, K+1]$. To cover in three intervals, the second interval thrown out must overlap the first or be contiguous to it: otherwise $\mathcal{C}$ will be two disjoint intervals and will have diameter more than $K$, whence one more step will not suffice to cover it. Again we may consider only the case where the second interval is contiguous to the right of the first and then double to count the case where the second is contiguous to the left of the first. The value of $j_2$ cannot be $j$ or less, since this would leave $\mathcal{C}$ with cardinality greater than $K+1$, which is too large a set to cover in one additional step. Thus, before doubling, the allowable range for $j_2$ is $[j+1, K+2]$. The corresponding range for $j_3$ is $[N-K, j_2+K+1]$. Equation (5.15) now becomes

$$
\begin{aligned}
(5.17) \quad \frac{\mathbb{P}(\mathcal{H}(3))}{2} &= \frac{1}{K+2} \frac{1}{N-1+(K+1)(K+1+j)} \\
&\quad \times \sum_{j_2=j+1}^{K+2} \frac{j_2-j}{N-2+(K+1)(K+2-j_2+j)} \\
&= \frac{1+o(1)}{K^3(1+j/K)} \sum_{s=1}^{K+2-j} \frac{1-j/K-s/K}{j+1+j/K+s},
\end{aligned}
$$

which is bounded when $j/K \in [\varepsilon, 1/2]$ and as $t := j/K \to 0^+$ due to

$$
\sum_{s=1}^{K+2-j} \frac{1-j/K-s/K}{j+1+j/K+s} \le \sum_{s=1}^{K+2-j} \frac{1}{j+1+j/K+s}
$$

$$
\approx \log(K(1-t)+2+Kt) - \log(Kt) = \log\left(\frac{K+2}{Kt}\right),
$$

by

$$
\frac{1+o(1)}{K^3} \log\left(\frac{1}{t}\right).
$$

Doubling yields, as a lower bound, the expression in the second line of Table 1 for $r=3$; for the upper bound, it remains to get an upper bound on $\mathbb{P}(\mathcal{H}^*(4))$.

We must sum this time over two types of sequences $(1, j_2, j_3)$. The first are those with $j+2 \notin [j+1, K+2]$; these do not appear in $\mathcal{H}(3)$ because it is not possible to cover $[N]$ in three intervals starting this way. The second are sequences where $j_2 \in [j+1, K+2]$ but $(1, j_2, j_3) \notin S(3)$; these do not appear in $\mathcal{H}(3)$ because the third interval did not complete the cover of $[N]$, where a different choice of $j_3$ could have completed the cover. Analyzing the second of these two types repeats the analysis from the last paragraph of Case 4. That is, allowing these values of $j_3$ replaces $1-j/K-s/K$ in the numerator of (5.17) by $1+j/K+s/K$, which does not affect the leading term when $j = o(K)$ and otherwise multiplies by a bounded factor, which we absorb into the definition of $f_3$.

The first of the two types of sequences splits into subtypes: $-j \le j_2 \le j$ (in which case you do not cover enough new ground to be able to complete coverage in three steps) or $K + 3 \leq j_2 \leq K + 1 + j$ [in which case the set $\mathsf{missed}(2, \mathbf{j})$ splits into two intervals and cannot be covered by one more interval]. For the first subtype, $M(3, \mathbf{j})$ is always at least $K$, so the sum over sequences of this subtype is $O(K^{-3})$. For the second subtype, $M(3, \mathbf{j}) = j$. For each $t$ there is exactly one value of $j_2$ for which $\mathsf{missed}(2, \mathbf{j})$ is composed of disjoint intervals sizes $t$ and $j - t$ in that order. Given that this occurs for some $t$, one may reason as in (5.3) to see that $Q(\mathbf{j}) \leq 2/(t(j - t))$. Thus the total probability of the second subtype is bounded above by

$$CK^{-3} \sum_{t=1}^{j-1} \frac{2}{t(j-t)} = O\left(K^{-3} \frac{\log j}{j}\right)$$

and since this is negligible compared to $K^{-3} \log(K/j)$, the proof in Case 3 is complete.

**6. The fat tail when $K = 1$.** In this section, we prove Theorem 2.6. For convenience we add a variable $Y_0$ to get an i.i.d. collection $\mathcal{C} := \{Y_0, Y_j, Y_{j,i} : 1 \leq j \leq N, 0 \leq i \leq 1\}$ and define the event $\widetilde{\mathcal{H}}_0$ to hold when $Y_0 \vee Y_1 \geq Y_{1,0} \vee Y_{1,1}$. Letting $\mathcal{H}^* = \widetilde{\mathcal{H}}_0 \cap \bigcap_{j=1}^{N} \mathcal{H}'_j$, it is evident that

$$\mathbb{P}(\mathcal{H}^*) \geq p_{\mathrm{fat}}(N + 1, 1)$$

by monotonicity of probability, and from Harris' (positive association) inequality we see that

$$p_{\mathrm{fat}}(N + 2, 1) \geq \mathbb{P}(\mathcal{H}'_{N+2} \cap \mathcal{H}'_{N+1} \cap \cdots \cap \mathcal{H}'_2) \cdot \mathbb{P}(\mathcal{H}'_1) = \mathbb{P}(\mathcal{H}^*) \cdot \mathbb{P}(\mathcal{H}'_1).$$

Since $\mathbb{P}(\mathcal{H}'_1) = c > 0$ independently of $N$, it suffices to prove Theorem 2.6 for $p_N := \mathbb{P}(\mathcal{H}^*)$ in places of $p_{\mathrm{fat}}(N, 1)$.

Having sliced open the circle, it is possible to derive a recursion for $p_N$. Observe that the order of the variables in $\mathcal{C}$, namely $\{Y_j, Y_{j,i}, Y_0 : 1 \leq j \leq N, 0 \leq i \leq 1\}$, is uniform among the $(3N + 1)!$ permutations, and that the permutation determines whether $\mathcal{H}^*$ has occurred. For $\mathcal{H}^*$ to occur, it is necessary that the maximum $M$ of variables in $\mathcal{C}$ be $Y_j$ for some $j$. Thus

$$p_N = \sum_{j=0}^{N} \frac{1}{3N + 1} \mathbb{P}(\mathcal{H}^* | Y_j = M).$$

These conditional probabilities may be evaluated recursively. If $Y_0 = M$, then further information about $Y_{1,0}$ and $Y_{1,1}$ is irrelevant and the ordering of the remaining $3N - 2$ variables is uniform, leading to

$$\mathbb{P}(\mathcal{H}^* | Y_0 = M) = p_{N-1}.$$

To ensure this holds for $N = 1$, we set $p_0 := 1$. Similarly,

$$\mathbb{P}(\mathcal{H}^* | Y_N = M) = p_{N-1}.$$

Now suppose $N \geq 2$ and $Y_j = M$ for some $1 < j < N-1$. Then $\mathcal{H}'_j$ and $\mathcal{H}'_{j+1}$ are known to occur. Removing from consideration the variables $Y_j$, $Y_{j,i}$ and $Y_{j+1,i}$ for $i = 0, 1$, the remaining variables are broken into two subsets of size $3(j-1)+1$ and $3(N-j-1)+1$; the ordering on the union of these is still jointly uniform, leading to

$$\mathbb{P}(\mathcal{H}^* | Y_j = M) = p_{j-1} p_{N-j-1}.$$

This equation is readily verified for $N \geq 2$ and $j = 1$ or $j = N-1$ as well. Putting these together gives the recursion

$$\begin{aligned}
(6.1) \qquad p_N &= \delta_{0,N} + \frac{1}{3N+1}\left(2p_{N-1} + \sum_{j=1}^{N-1} p_{j-1} p_{N-j-1}\right) \\
&= \frac{1}{3N+1}\left(2p_{N-1} + \sum_{j=2}^{N} p_{j-2} p_{N-j}\right),
\end{aligned}$$

which holds for all $N$ due to the inclusion of the delta function.

Let $f(z) := \sum_{N=0}^{\infty} p_N z^N$. Since we know (by submultiplicativity) that $\log p_N / N \to \log(\lambda)$ for some $\lambda \in (0, 1)$, the radius of convergence for the power series defining $f$ above will be $1/\lambda$. The generating function for $(3N+1)p_N$ is equal to $f + 3zf'$. The generating function for $\delta_{0,N}$ is 1, the generating function for $2p_{N-1}$ is $2zf$ and the generating function for $\sum_{j=2}^{N} p_{j-2} p_{N-j}$ is $z^2 f^2$. Equation (6.1) then becomes a Riccati equation:

$$(6.2) \qquad f + 3zf' = 1 + 2zf + z^2 f^2.$$

From the derivation it is apparent that this functional equation has a unique formal power series solution, $f$, and since $|p_N| \leq 1$ for all $N$, the series represents a function, also denoted $f$, that is analytic in a neighborhood of the origin. Only one locally analytic function can satisfy (6.2). To see this, write $g(z) = zf(z^3)$ so that $g' = 1 + 2z^2 g + z^4 g^2 := F(z, g)$ with boundary value $g(0) = 0$. Since $F$ is bounded and Lipschitz in a neighborhood of the origin, Gronwall's lemma ([5] or implicit in the classical uniqueness result [1], Theorem 2.2) says there is at most one such $g$ in the set of functions differentiable near 0.

Thus $f$ is the unique locally analytic solution to (6.2), whence we may use Maple's ordinary differential equation solver to find solutions to (6.2) and be rigorously assured that any such solution we can verify by differentiation must equal $f$. One finds that for any real constant $A$, there is a solution $f_A$ which is a ratio of Bessel functions. Its numerator is equal to

$$\mathbf{num} := (A \operatorname{BesselI}(-\tfrac{1}{3}, \tfrac{2}{3}\sqrt{2z}) + \operatorname{BesselK}(\tfrac{1}{3}, \tfrac{2}{3}\sqrt{2z}))$$

and its denominator is equal to

$$\mathbf{den} := \sqrt{z}(-A\sqrt{2}\ \mathrm{BesselI}(\tfrac{2}{3}, \tfrac{2}{3}\sqrt{2z}) + A\sqrt{z}\ \mathrm{BesselI}(-\tfrac{1}{3}, \tfrac{2}{3}\sqrt{2z})$$
$$+ \sqrt{2}\ \mathrm{BesselK}(\tfrac{2}{3}, \tfrac{2}{3}\sqrt{2z}) + \sqrt{z}\ \mathrm{BesselK}(\tfrac{1}{3}, \tfrac{2}{3}\sqrt{2z}));$$

here BesselI and BesselK denote modified Bessel functions of the first and second kinds, respectively. It is not yet clear whether one of these solutions is $f$.

As a fractional power series, $f_A$ has a leading term of $z^{-1/3}$, so certainly if $f_A = f$, then $A$ must be chosen to make this term vanish. Solving for $A$ yields $A = -\pi\sqrt{3}/3$, and plugging this into the expressions for $\mathbf{num}$ and $\mathbf{den}$ leads to a function with a power series, a priori fractional, beginning with $1 + z/2 + \cdots$. The series converges in a neighborhood of the origin, so it defines a function that is $1 + O(z)$ near $z = 0$. Any function that is $1 + O(z)$ near the origin and satisfies the differential equation (6.2) must be analytic in a neighborhood of the origin. We have therefore found the function $f$.

Since $f$ has positive coefficients, its minimal modulus singularities lie on the positive real axis. Its functional form dictates that $f$ has positive real singularities precisely at the zeros of $\mathbf{den}$. We may approximate these as closely as we wish. Maple's numeric solver gives $z_0 := 1.803034611\ldots$ (the constant is not recognized by Plouffe's inverse symbolic calculator). Thus

$$\frac{\log p_N}{N} \to -\log z_0 = -0.58947114\ldots,$$

which finishes the proof of Theorem 2.6.

**Acknowledgment.** We are grateful to Yuval Peres for useful suggestions concerning the analysis of the normal case.

## REFERENCES

[1] CODDINGTON, E. and LEVINSON, N. (1955). *Theory of Ordinary Differential Equations.* McGraw-Hill, New York. MR69338

[2] DURRETT, R. and LIMIC, V. (2003). Rigorous results for the NK model. *Ann. Probab.* **31** 1713–1753. MR2016598

[3] EVANS, S. and STEINSALTZ, D. (2002). Estimating some features of NK fitness landscapes. *Ann. Appl. Probab.* **12** 1299–1321. MR1936594

[4] FELLER, W. (1971). *An Introduction to Probability Theory and Its Applications* **2**, 2nd ed. Wiley, New York.

[5] HIRSCH, M. and SMALE, S. (1974). *Differential Equations, Dynamical Systems, and Linear Algebra.* Academic Press, New York. MR486784

[6] KAUFFMAN, S. (1993). *The Origins of Order.* Oxford Univ. Press.

[7] KAUFFMAN, S. and LEVIN, S. (1987). Towards a general theory of adaptive walks on rugged landscapes. *J. Theoret. Biol.* **128** 11–45. MR907587

[8] KNIGHT, F. (1981). *Essentials of Brownian Motion and Diffusion.* Amer. Math. Soc. Providence, RI. MR613983

[9] Revuz, D. and Yor, M. (1994). *Continuous Martingales and Brownian Motion*, 2nd ed. Springer, New York. MR1303781

[10] Stanley, R. P. (1986). *Enumerative Combinatorics*, I. Wadsworth and Brooks/Cole, Belmont, CA. MR847717

[11] Weinberger, E. (1991). Local properties of Kauffman's NK model: A tunably rugged energy landscape. *Phys. Rev. A* **44** 6399–6413.

Department of Mathematics
University of British Columbia
121-1984 Mathematics Road
Vancouver, British Columbia
Canada V6T 1Z2
e-mail: limic@math.ubc.ca

Department of Mathematics
University of Pennsylvania
209 S. 33rd Street
Philadelphia, Pennsylvania 19104-6395
USA
e-mail: pemantle@math.upenn.edu